\documentclass{svproc}
\usepackage{url}

\usepackage{graphicx}
\usepackage{multicol}
\usepackage{footmisc}
\usepackage{subfigure}
\usepackage{amsfonts,amsmath,amssymb}
\usepackage[latin1]{inputenc}
\usepackage{dsfont}
\usepackage{enumerate}
\usepackage{xcolor}
\usepackage[all]{xy}
\def\notshow#1\notshowend{} %
\usepackage{hyphenat}
\usepackage{microtype}

\usepackage{graphicx}

\usepackage{tikz}
\usepackage{tikz-3dplot}

\usepackage{multirow}
\usepackage{array}

\usepackage{amssymb}
\usepackage{amsfonts}
\usepackage{mathrsfs}
\usepackage{hyperref}
\usepackage{mathtools}

\usepackage[normalem]{ulem} 

\def\bb#1\eb{\textcolor{blue}{#1}} 
\def\br#1\er{\textcolor{red}{#1}} %
\def\bm#1\em{\textcolor{purple}{#1}} %

\DeclareMathOperator{\atantwo}{atan2}

\newtheorem{thm}{Theorem}
\newtheorem{prop}{Proposition}

\newtheorem{defi}{Definition}

\newtheorem{exe}{Example}

\newtheorem{rem}{Remark}

\newcommand{\ben}{\begin{enumerate}}
\newcommand{\een}{\end{enumerate}}
\newcommand{\bit}{\begin{itemize}}
\newcommand{\eit}{\end{itemize}}
\newcommand{\edoc}{\end{document}}

\begin{document}
\mainmatter              
\title{On the application of Lorentz-Finsler geometry to model wave propagation\thanks{This is a preprint of the following chapter: E.~Pend\'{a}s-Recondo, On the Application of Lorentz-Finsler Geometry to Model Wave Propagation, published in Progress in Lorentzian Geometry, Springer Proceedings in Mathematics \& Statistics, vol. 512, edited by W.~Barrera, J.~Herrera, J.~P.~Navarrete, M.~Navarro, O.~Palmas and D.~A.~Solis, 2025, Springer, reproduced with permission of Springer Nature Switzerland AG. The final authenticated version is available online at: \url{https://doi.org/10.1007/978-3-031-99212-4}.}}
\titlerunning{On the application of Lorentz-Finsler geometry to model wave propagation}  
%
\author{Enrique Pend\'{a}s-Recondo}
\authorrunning{E. Pend\'{a}s-Recondo} 
%
%
\institute{Departamento de Matem\'{a}ticas, Universidad de Murcia, Campus de Espinardo, 30100 Espinardo, Murcia, Spain. \email{e.pendasrecondo@um.es}.}

\maketitle              

\begin{abstract}
The recent increasing interest in the study of Lorentz-Finsler geometry has led to several applications to model real-world physical phenomena. Our purpose is to provide a simple, step-by-step review on how to build and implement such a geometric model to describe the propagation of a classical wave satisfying Fermat's and Huygens' principles in an anisotropic and rheonomic (time-dependent) medium. The model is based on identifying the individual wave trajectories as lightlike pregeodesics of a specific Lorentz-Finsler metric, which obey a simple ODE system and can therefore be easily computed in real time.

\end{abstract}

\newpage

\section{Introduction}
Finsler geometry has proven to be a very powerful tool to model, in a very natural way, physical problems that exhibit an intrinsic anisotropy. The recent growing interest in the study of these models began with Zermelo's navigation problem, first proposed in \cite{zermelo1931}, but then revisited and solved in \cite{bao2004,shen2003} from a geometric point of view. In this problem, one seeks the fastest trajectory between two points for a moving object in the presence of a current. In the time-independent case, it turns out that the set of maximum velocity vectors of the object (taking into account the current) gives rise to a Finsler metric $ F $---essentially, a smooth distribution of non-symmetric norms---that measures the traveltime: the $ F $-length of a curve coincides with the time spent by the moving object following this path at maximum speed. The solution to Zermelo's problem, then, can be characterized as a pregeodesic curve of $ F $ (i.e., a critical point of the traveltime).\footnote{A complete study of this problem via Finsler geometry can be found in \cite{caponio2014}.}

This way of interpreting Zermelo's problem provides an interesting connection with the classical Fermat's principle, which states that the path taken by a wave ray between two points is precisely the one that makes critical the traveltime. So, if we consider a classical wave that spreads in an anisotropic medium, each individual wave trajectory can also be interpreted as a pregeodesic of a certain Finsler metric. This has led to several applications to real-world physical phenomena---especially when there is an anisotropy involved---such as the modeling of sound waves (see \cite{gibbons2010,gibbons2011}), seismic waves (see \cite{antonelli2003,yajima2009}) or wildfire propagation (see \cite{dehkordi2022,markvorsen2016}).

However, one last step to complete these models is to allow the speed of the wave to vary over time, generating a time-dependent Finsler metric. This greatly complicates the description of the problem from a purely Finslerian viewpoint due to the fact that, now, pregeodesic curves on the space cannot be defined as usual. The way to overcome this difficulty is to add explicitly the time dimension to the space and define a Lorentz-Finsler metric whose lightcone (the vectors with vanishing norm) represents the infinitesimal propagation of the wave through the resulting spacetime---in the spirit of the way general relativity models the (isotropic) propagation of light using Lorentzian geometry. Then, wave trajectories can be characterized as lightlike pregeodesics of this Lorentz-Finsler metric. This provides a unified geometric framework to model wave propagation under general time-, position- and direction-dependent conditions, first developed in \cite{pendas2021} and then applied specifically to wildfires in \cite{pendas2023a,pendas2024}.

The aim of this work is to provide a simple, step-by-step review on how to build this model and how it can be implemented in practice---which may serve as an introduction to more technical references. We begin in \S~\ref{sec:setting} by describing the different ingredients of the model: the measure of the distance (\S~\ref{subsec:distance}), the speed profile of the wave (\S~\ref{subsec:speed}) and how it gives rise to a Finsler metric (\S~\ref{subsec:F}), the construction of the spacetime and the Lorentz-Finsler metric in the time-dependent case (\S~\ref{subsec:L}), the measure of the traveltime (\S~\ref{subsec:traveltime}), the initial wavefront (\S~\ref{subsec:initial_wavefront}) and the deduction of the geodesic equations (\S~\ref{subsec:geodesics}). Then, in \S~\ref{sec:model} we explicitly construct the model, identifying the wave trajectories (\S~\ref{subsec:wave_trajectories}) and the wavefront (\S~\ref{subsec:wavefront}), to finally describe its implementation in practice (\S~\ref{subsec:implementation}). We end in \S~\ref{subsec:advantages} by briefly discussing the main advantages that this model provides.

\section{The setting}
\label{sec:setting}
The main goal is to build a model that effectively solves the following problem:
\begin{quote}
Knowing the speed of a classical wave at each point, time and direction, can we compute the entire wave propagation from a given initial wavefront?
\end{quote}
Let us proceed step by step, constructing the setting needed to answer this question.

\subsection{Distance}
\label{subsec:distance}
The first thing we need to specify is the {\em space} over which the wave spreads. This space will be assumed to be a smooth manifold of dimension $ n $. Since we are interested on the applications of the model and explicit computations will be needed, we will work in coordinates in $ N \subset \mathds{R}^n $. Namely, a point in the space will be written as $ x = (x^1,\ldots,x^n) \in N \subset \mathds{R}^n $ and a tangent vector as $ v = (v^1,\ldots,v^n) \in T_xN \equiv \mathds{R}^n $. In real-world applications, this space $ N $ is always endowed with a Riemannian metric $ h $ (i.e., a scalar product $ h_x(\cdot,\cdot) $ on each tangent space $ T_xN $), which measures actual distances on $ N $. This means that if $ \sigma: [a,b] \rightarrow N $ is a curve on $ N $, the distance we travel when following $ \sigma $ is
\begin{equation}
\label{eq:length}
L_h[\sigma] \coloneqq \int_a^b \sqrt{h_{\sigma(s)}(\dot{\sigma}(s),\dot{\sigma}(s))} ds = \int_a^b ||\dot{\sigma}(s)||_{h_{\sigma(s)}} ds,
\end{equation}
where $ || v ||_{h_x} \coloneqq \sqrt{h_x(v,v)} $ denotes the norm defined by $ h $. This computation is independent of the parametrization of $ \sigma $.

\begin{exe}
\label{ex:metrics}
The most common situations when dealing with wave propagation are the following:
\begin{enumerate}
\item[(i)] The wave spreads directly over $ N \subset \mathds{R}^n $ (usually $ n=2,3 $) and the distance is measured by the natural Euclidean metric:
\begin{equation*}
h(v,u) = v^1 u^1 + \ldots + v^n u^n,
\end{equation*}
which is independent of the position $ x \in N $.
\item[(ii)] The wave spreads over a surface (in particular, a graph) of $ \mathds{R}^3 $ given by a height function $ z: N \subset \mathds{R}^2 \rightarrow \mathds{R}^3 $. Namely, each coordinate $ x = (x^1,x^2) \in N $ (the projection on the horizontal plane) univocally corresponds to a point $ (x^1,x^2,z(x^1,x^2)) \in \mathds{R}^3 $ of the actual surface where the propagation is taking place (see e.g. \cite{pendas2023a,markvorsen2016}). In this case, the distance is measured by the Riemmanian metric induced on the surface from the natural Euclidean metric on $ \mathds{R}^3 $:
\begin{equation}
\label{eq:metric2}
h_x(v,u) =
\begin{pmatrix}
v^1 & v^2
\end{pmatrix}
\begin{pmatrix}
1 + \left( \frac{\partial z}{\partial x^1} \right)^2 & \frac{\partial z}{\partial x^1} \frac{\partial z}{\partial x^2} \\
\frac{\partial z}{\partial x^2} \frac{\partial z}{\partial x^1} & 1 + \left( \frac{\partial z}{\partial x^2} \right)^2 \\
\end{pmatrix}
\begin{pmatrix}
u^1 \\
u^2 \\
\end{pmatrix}.
\end{equation}
Note that this Riemannian metric is position-dependent.
\end{enumerate}
\end{exe}

If we consider the set $ \mathcal{N}_{x,y} $ of all (piecewise smooth) curves $ \sigma: [a,b] \rightarrow N $ going from $ \sigma(a) = x $ to $ \sigma(b) = y $, the one that minimizes the distance between $ x $ and $ y $ must be a critical point (i.e., it satisfies the usual Euler-Lagrange equations) of the $ h $-{\em length functional} $ L_h $, given by \eqref{eq:length}. Individual wave rays, however, do not minimize (or more generally, make critical) the distance, but the {\em traveltime}, according to the classical {\em Fermat's principle}:
\begin{quote}
The path taken by a wave ray between two given points is the one that makes critical the  traveltime.
\end{quote}
Let us now identify these curves.

\subsection{Speed profile}
\label{subsec:speed}
The initial hypotheses assume that we know the {\em speed profile} of the wave, which can be, in the most general case, rheonomic (time-dependent), inhomogeneous (position-dependent) and anisotropic (direction-dependent). This means that there exists a positive function $ V(t,x,v) $ representing the speed of the wave at each (absolute) time $ t \in \mathds{R} $,\footnote{Here we will not consider relativistic effects, so $ t $ is absolute in the sense that it is the time measured by every observer.} each point $ x \in N $ and each direction $ v \in T_xN $. Obviously, this speed should not depend on the length of $ v $, only on the direction it points, so $ V $ must be positive homogeneous of degree $ 0 $ on the direction:\footnote{It will be assumed that the wave can spread in all directions. The case of a ``strong wind'' effect that makes the wave unable to propagate in some directions is studied in detail in \cite[\S~6]{pendas2021}.}
\begin{equation*}
V(t,x,\lambda v) = V(t,x,v), \quad \forall \lambda > 0.
\end{equation*}
For example, if $ n = 2 $ we can express the speed profile as $ V(t,x,\theta) $, where $ \theta \in [0,2\pi) $ denotes the angle in the usual polar coordinates (see Example~\ref{ex:finsler_metrics}(i) below).

Among all the curves in $ \mathcal{N}_{x,y} $, assume that $ \sigma: [a,b] \rightarrow N $ is the path taken by an individual wave ray. In general, $ \sigma(s) $ can be arbitrarily parametrized, but we can always consider a time reparametrization $ s(t) $, obtaining a new curve $ \sigma \circ s(t) = \sigma(s(t)) $ whose parameter is the actual time coordinate $ t $, i.e. $ \sigma \circ s $ marks the actual position of the wave ray over time. Physically, the speed is the magnitude (in our case, the $ h $-norm) of the velocity vector, which is in turn the derivative of the position with respect to the time. Therefore:
\begin{equation*}
V(t(s),\sigma(s),\dot{\sigma}(s)) = \left\lvert\left\lvert \frac{d(\sigma \circ s)}{dt}(t(s)) \right\rvert\right\rvert_{h_{\sigma(s)}} = \dot{s}(t(s)) ||\dot{\sigma}(s)||_{h_{\sigma(s)}},
\end{equation*}
and noting that $ \dot{s}(t(s)) = \frac{1}{\dot{t}(s)} $, we obtain
\begin{equation}
\label{eq:time_der}
\dot{t}(s) = \frac{||\dot{\sigma}(s)||_{h_{\sigma(s)}}}{V(t(s),\sigma(s),\dot{\sigma}(s))},
\end{equation}
which is a necessary condition for $ \sigma \in \mathcal{N}_{x,y} $ to represent a wave path.

\subsection{Finsler metrics}
\label{subsec:F}
Let us take a closer look at the right-hand side of \eqref{eq:time_der}. It defines a function
\begin{equation}
\label{eq:finsler}
F_{t,x}(v) \coloneqq \frac{||v||_{h_x}}{V(t,x,v)}, \quad \forall v \in T_xN
\end{equation}
(with $ F_{t,x}(0) \coloneqq 0 $), which keeps some of the basic properties of a norm, with the particularity that it is time-dependent:
\begin{itemize}
\item $ F $ is positive: $ F_{t,x}(v) > 0 $, for all $ t \in \mathds{R} $, $ x \in N $ and $ v \in T_xN \setminus \{0\} $.
\item $ F $ is positive homogeneous of degree 1: $ F_{t,x}(\lambda v) = \lambda F_{t,x}(v) $, for all $ \lambda > 0 $.
\item $ F $ is smooth away from the zero vector, as long as $ V $ is smooth.
\end{itemize}
Also, related to this function, for each fixed time and position $ (t,x) \in \mathds{R} \times N $:
\begin{itemize}
\item We define the {\em indicatrix} $ \Sigma_{t,x} $ of $ F_{t,x} $ as the set of $ F_{t,x} $-unit vectors:
\begin{equation*}
\Sigma_{t,x} \coloneqq \{v \in T_xN: F_{t,x}(v) = 1\},
\end{equation*}
which is a hypersurface of $ T_x N $ diffeomorphic to a sphere and enclosing the zero vector (e.g., when $ n = 2 $, the indicatrix is a simple closed curve enclosing the origin).
\item We define the {\em fundamental tensor} of $ F_{t,x} $ in the direction $ v \in T_xN \setminus \{0\} $, denoted $ g_v^{F_{t,x}} $, as the symmetric bilinear form given by the Hessian of $ \frac{1}{2}F^2_{t,x} $ at $ v $:\footnote{This is the way one obtains a metric out of a norm. For instance, the Hessian of $ \frac{1}{2}||\cdot||_h^2 $ provides the Riemannian metric $ h $.}
\begin{equation*}
g^{F_{t,x}}_v(u,w) \coloneqq \frac{1}{2}\textup{Hess}(F_{t,x}(v)^2)(u,w) = \frac{1}{2} \left. \frac{\partial^2}{\partial \delta \partial \eta} F_{t,x}(v + \delta u + \eta w)^2 \right\rvert_{\delta=\eta=0}.
\end{equation*}
In matrix form, using coordinates:
\begin{equation}
\label{eq:gF_coordinates}
g^{F_{t,x}}_v(u,w) = \frac{1}{2}
\begin{pmatrix}
u^1 & \cdots & u^n
\end{pmatrix}
\begin{pmatrix}
\frac{\partial^2 F_{t,x}^2}{\partial(v^1)^2}(v) & \cdots & \frac{\partial^2 F_{t,x}^2}{\partial v^1 \partial v^n}(v) \\
\vdots & \ddots & \vdots \\
\frac{\partial^2 F_{t,x}^2}{\partial v^n \partial v^1}(v) & \cdots & \frac{\partial^2 F_{t,x}^2}{\partial(v^n)^2}(v)
\end{pmatrix}
\begin{pmatrix}
w^1 \\
\vdots \\
w^n
\end{pmatrix},
\end{equation}
and notice that $ g^{F_{t,x}}_v(v,v) = F_{t,x}(v)^2 $.
\end{itemize}

These two objects are related to each other by the following result (see e.g. \cite[Proposition~2.3(v)]{javaloyes2014}).
\begin{prop}
\label{prop:finsler}
For any fixed $ (t,x) \in \mathds{R} \times N $, the following statements are equivalent:
\begin{itemize}
\item The fundamental tensor $ g^{F_{t,x}}_v(\cdot,\cdot) $ is a positive definite scalar product on $ T_xN $, for all $ v \in T_xN \setminus \{0\} $.
\item The indicatrix $ \Sigma_{t,x} \subset T_xN $ is strongly convex, i.e. using the natural Euclidean metric on $ T_xN \equiv \mathds{R}^n $, the second fundamental form of $ \Sigma_{t,x} $ with respect to the inner normal vector is positive definite everywhere.
\item $ \Sigma_{t,x} $ has positive sectional curvature everywhere with respect to the natural Euclidean metric on $ T_xN \equiv \mathds{R}^n $.
\end{itemize}
If this holds for every $ (t,x) \in \mathds{R} \times N $, then $ F $ is said to be a {\em (time-dependent) Finsler metric} on $ N $.\footnote{See \cite{javaloyes2014} for the formal definition and properties of a Finsler metric.}
\end{prop}

When $ F $ is a Finsler metric, it satisfies the strict triangle inequality, which is the last basic property usually required to norms:\footnote{In fact, this inequality holds as long as $ \Sigma_{t,x} $ is strictly convex (see \cite[Proposition~2.3(iv)]{javaloyes2014}).}
\begin{equation*}
F_{t,x}(v+u) \leq F_{t,x}(v) + F_{t,x}(u), \quad \forall v,u \in T_xN \setminus \{0\},
\end{equation*}
with equality if and only if $ v = \lambda u $ for some $ \lambda > 0 $. This way, fixing $ (t,x) \in \mathds{R} \times N $, we can think of $ F_{t,x} $ as a non-symmetric norm (in general $ F_{t,x}(v) \not= F_{t,x}(-v) $), the fundamental tensor $ g_v^{F_{t,x}} $ as its associated scalar product for each direction $ v $, and the indicatrix $ \Sigma_{t,x} $ represents the infinitesimal propagation of the wave: $ \Sigma_{t,x} $ is the wavefront one would obtain after one time unit if the wave departed from $ (t,x) $ under constant conditions (i.e., if $ F $ were time- and position-independent).

\begin{rem}
In the following, we will always assume that the speed profile of the wave is such that the statements in Proposition~\ref{prop:finsler} hold, i.e. $ F $ in \eqref{eq:finsler} is a (time-dependent) Finsler metric on $ N $. This is not very restrictive when modeling wave propagation, as it only requires the infinitesimal propagation to be strongly convex (see the discussion in \cite[Appendix~A.2]{pendas2023a}). The complete propagation, on the other hand, can be non-convex.
\end{rem}

\begin{exe}
\label{ex:finsler_metrics}
Using the Riemannian metrics of Example~\ref{ex:metrics}, some typical speed profiles and Finsler metrics are the following:
\begin{itemize}
\item[(i)] On $ N = \mathds{R}^2 $, consider the speed profile
\begin{equation*}
V(t,x,\theta) = \frac{a(t,x)(1-\varepsilon(t,x)^2)}{1-\varepsilon(t,x) \cos(\theta-\phi(t,x))},
\end{equation*}
for some smooth positive functions $ a(t,x) $, $ \varepsilon(t,x) $ and $ \phi(t,x) $. This provides the (time-dependent) Finsler metric
\begin{equation*}
F_{t,x}(v) = \frac{\sqrt{(v^1)^2+(v^2)^2}}{a(t,x)(1-\varepsilon(t,x)^2)} \left[ 1-\varepsilon(t,x) \cos\left(\atantwo(v^2,v^1)-\phi(t,x)\right) \right],
\end{equation*}
whose indicatrix $ \Sigma_{t,x} \subset \mathds{R}^2 $ is an ellipse of eccentricity $ \varepsilon(t,x) $ centered at one of their foci, with semi-major axis $ a(t,x) $ oriented in the direction $ \phi(t,x) $. An elliptical propagation of this type is typically used to model the spread of wildfires in the presence of wind (see \cite{finney1998,pendas2023a}).\footnote{Another example for the choice of $ V(t,x,\theta) $ in $ \mathds{R}^2 $ is Gielis superformula, which features a wide range of different shapes (see \cite{pendas2024}).}

\item[(ii)] On the projection $ N \subset \mathds{R}^2 $ of a surface of $ \mathds{R}^3 $ given by the height function $ z(x^1,x^2) $, consider the speed profile
\begin{equation*}
V(t,x,v) = b(t,x) \pm c(t,x) \frac{v^1 \frac{\partial z}{\partial x^1} + v^2 \frac{\partial z}{\partial x^2}}{\sqrt{h_x(v,v)}},
\end{equation*}
for some positive functions $ b(t,x) $ and $ c(t,x) $, and the induced Riemannian metric $ h $ defined in \eqref{eq:metric2}. This provides the (time-dependent) Finsler metric
\begin{equation*}
F_{t,x}(v) = \frac{h_x(v,v)}{b(t,x) \sqrt{h_x(v,v)} \pm c(t,x)\left( v^1 \frac{\partial z}{\partial x^1} + v^2 \frac{\partial z}{\partial x^2} \right)},
\end{equation*}
which is a Matsumoto metric (see e.g. \cite{shimada2005}), used to model wave propagation that takes place on a surface and is affected by the slope. Choosing the $ + $ sign means that the propagation is faster upwards (as is the case for wildfires; see \cite[\S~3.1]{pendas2023a}), whereas the $ - $ sign means that it is faster downwards.
\end{itemize}
\end{exe}

Geodesic curves will play a key role in this model, as we will see later. When $ F $ is time-independent, we say that it is a {\em proper Finsler metric} on $ N $ and, similar to the case with Riemannian metrics, we can define the {\em geodesics} of $ F $ as the critical points of the {\em $ F $-energy functional}
\begin{equation*}
E_F[\sigma] \coloneqq \int_a^b F_{\sigma(s)}(\dot{\sigma}(s))^2 ds,
\end{equation*}
among all (piecewise smooth) curves $ \sigma: [a,b] \rightarrow N $ with fixed endpoints (see e.g. \cite{caponio2011}). However, these critical points are preserved only under positive affine reparametrizations of $ \sigma $,\footnote{Specifically, $ \sigma(s) $ is a geodesic of $ F $ if and only if so is $ \sigma \circ s(r) $, with $ s(r) = \lambda r + c $ for some $ \lambda > 0 $ and $ c \in \mathds{R} $.} i.e. geodesics require a specific parametrization---essentially, they represent trajectories without acceleration---whereas here we are simply interested on the path they follow. In order to remove this dependence on the parametrization, we define the {\em pregeodesics} of $ F $ as the curves that can be reparametrized as geodesics (keeping the orientation). It turns out that pregeodesics are the critical points of the {\em $ F $-length functional} (see \cite[\S~5]{bao2000})
\begin{equation}
\label{eq:F-length}
L_F[\sigma] \coloneqq \int_a^b F_{\sigma(s)}(\dot{\sigma}(s)) ds,
\end{equation}
which is indeed independent of reparametrizations of $ \sigma $, thanks to the fact that $ F $ is positive homogeneous of degree 1.

In the general time-dependent case, however, $ E_F $ and $ L_F $ not only depend on $ \sigma(s) $ but also on the external time function $ t(s) $. Therefore, (pre-)geodesic curves cannot be defined as usual---the same curve $ \sigma $ on $ N $ can have different $ F $-energy and $ F $-length in different contexts---and we need to generalize this notion.

\subsection{Lorentz-Finsler metrics}
\label{subsec:L}
To deal with the time dependence and solve the issue of defining geodesic curves in this case, time must be included as an intrinsic element of the manifold, rather than as an external function. This can be achieved by explicitly treating time as an additional dimension and considering the {\em spacetime} $ M \coloneqq \mathds{R} \times N $, where the natural projection (time coordinate) $ t: M \rightarrow \mathds{R} $ measures the absolute time. This way, points in the spacetime will be expressed as $ p = (t,x) \in M $, vectors as $ \hat{v} = (v^0,v) \in T_pM \equiv \mathds{R} \times T_xN $ and curves $ \gamma: [a,b] \rightarrow M $ as $ \gamma(s) = (t(s),\sigma(s)) $. The idea then is to mimic the way general relativity models the propagation of light, constructing a metric with Lorentzian signature $ (+,-,\ldots,-) $ where curves with vanishing norm represent trajectories of individual wave rays.

From \eqref{eq:time_der}, we know that if $ \gamma(s) = (t(s),\sigma(s)) $ is the path taken by a wave ray in the spacetime, then
\begin{equation}
\label{eq:time_der_F}
\dot{t}(s) = F_{\gamma(s)}(\dot{\sigma}(s)).
\end{equation}
Since $ \dot{t}(s) > 0 $ (time increases as the wave ray travels) and $ F $ is positive, this is equivalent to
\begin{equation*}
\dot{t}(s)^2 - F_{\gamma(s)}(\dot{\sigma}(s))^2 = 0,
\end{equation*}
so if we consider the function
\begin{equation}
\label{eq:def_L}
G_p(\hat{v}) \coloneqq (v^0)^2 - F_p(v)^2, \quad \forall \hat{v} = (v^0,v) \in T_pM,
\end{equation}
we can rewrite condition \eqref{eq:time_der_F} as
\begin{equation}
\label{eq:time_der_L}
G_{\gamma(s)}(\dot{\gamma}(s)) = 0.
\end{equation}

\begin{defi}
The function $ G $ defined in \eqref{eq:def_L}, being $ F $ a (time-dependent) Finsler metric, is called a {\em Lorentz-Finsler metric} with {\em lightcone}
\begin{equation*}
\mathcal{C}_p \coloneqq \{\hat{v} = (v^0,v) \in T_pM: G_p(\hat{v}) = 0 \text{ and } v^0 > 0 \}, \quad \forall p \in M.
\end{equation*}
Moreover, we say that a vector $ \hat{v} \in T_pM $ is {\em lightlike} if $ \hat{v} \in \mathcal{C}_p $ and, similarly, a curve $ \gamma(s) $ is {\em lightlike} if so is $ \dot{\gamma}(s) $ everywhere.\footnote{See \cite{javaloyes2020} for the formal definition and properties of Lorentz-Finsler metrics and Finsler spacetimes. We will assume that this spacetime is globally hyperbolic, which is not restrictive for modeling (see \cite[Remark~3.2 and Convention~3.3]{pendas2021}).}
\end{defi}

\begin{rem}
\label{rem:lightlike}
Observe that, for any $ v \in T_xN \setminus \{0\} $ and $ t \in \mathds{R} $, there exists a unique $ v^0 > 0 $ such that $ \hat{v} = (v^0,v) \in \mathcal{C}_{t,x} $. In fact, $ v \in \Sigma_{t,x} $ if and only if $ (1,v) \in \mathcal{C}_{t,x} $. 
\end{rem}

Therefore, from the previous observations, we can conclude that wave rays follow lightlike curves and the lightcone $ \mathcal{C}_p $ represents the infinitesimal propagation of the wave through the spacetime, much like $ \Sigma_{t,x} $ represents it through the space. Note also the difference between the properties of $ G $ and those of $ F $ in \S~\ref{subsec:F}:
\begin{itemize}
\item $ G $ is positive homogeneous of degree 2: $ G_p(\lambda \hat{v}) = \lambda^2 G_p(\hat{v}) $, for all $ \lambda > 0 $. The reason to use a 2-homogeneous function instead of a 1-homogeneous one such as $ F $ (and thus, the reason why we squared equation \eqref{eq:time_der_F} to obtain $ G $) is because 1-homogeneous functions fail to be smooth on lightlike vectors.
\item $ G $ is smooth away from vectors of the form $ \hat{v} = (\lambda,0) $, with $ \lambda \in \mathds{R} $ (unless $ F $ is Riemannian, in which case $ G $ is a smooth Lorentzian norm).
\item We define the {\em fundamental tensor} of $ G_p $ in the direction $ \hat{v} \in \mathcal{C}_p $, denoted $ g_{\hat{v}}^{G_p} $, as the Hessian of $ \frac{1}{2}G_p $ at $ \hat{v} $:
\begin{equation*}
\begin{split}
g^{G_p}_{\hat{v}}(\hat{u},\hat{w}) & \coloneqq \frac{1}{2}\textup{Hess}(G_p(\hat{v}))(\hat{u},\hat{w}) = \frac{1}{2} \left. \frac{\partial^2}{\partial \delta \partial \eta} G_p(\hat{v} + \delta \hat{u} + \eta \hat{w}) \right\rvert_{\delta=\eta=0} \\
& = u^0 w^0 - g^{F_p}_{v}(u,w).
\end{split}
\end{equation*}
In matrix form, using coordinates:
\begin{equation}
\label{eq:gL_coordinates}
\begin{split}
g^{G_p}_{\hat{v}}(\hat{u},\hat{w}) & = \frac{1}{2}
\begin{pmatrix}
u^0 & \cdots & u^n
\end{pmatrix}
\begin{pmatrix}
\frac{\partial^2 G_p}{\partial(v^0)^2}(\hat{v}) & \cdots & \frac{\partial^2 G_p}{\partial v^0 \partial v^n}(\hat{v}) \\
\vdots & \ddots & \vdots \\
\frac{\partial^2 G_p}{\partial v^n \partial v^0}(\hat{v}) & \cdots & \frac{\partial^2 G_p}{\partial(v^n)^2}(\hat{v}) 
\end{pmatrix}
\begin{pmatrix}
w^0 \\
\vdots \\
w^n
\end{pmatrix} \\
& =
\begin{pmatrix}
u^0 & \cdots & u^n
\end{pmatrix}
\begin{pmatrix}
1 & 0 & \ldots & 0 \\
0 & -\frac{1}{2}\frac{\partial^2 F^2_p}{\partial(v^1)^2}(v) & \cdots & -\frac{1}{2}\frac{\partial^2 F^2_p}{\partial v^1 \partial v^n}(v) \\
\vdots & \vdots & \ddots & \vdots \\
0 & -\frac{1}{2}\frac{\partial^2 F^2_p}{\partial v^n \partial v^1}(v) & \cdots & -\frac{1}{2}\frac{\partial^2 F^2_p}{\partial(v^n)^2}(v)
\end{pmatrix}
\begin{pmatrix}
w^0 \\
\vdots \\
w^n
\end{pmatrix}.
\end{split}
\end{equation}
Note that $ g^{G_p}_{\hat{v}}(\hat{v},\hat{v}) = G_p(\hat{v}) $ and $ g^{G_p}_{\hat{v}}(\cdot,\cdot) $ is a Lorentzian scalar product with signature $ (+,-,\ldots,-) $.
\end{itemize}

Now we can define the {\em geodesics} of $ G $ as the critical points of the {\em $ G $-energy functional}
\begin{equation}
\label{eq:L-energy}
E_G[\gamma] \coloneqq \int_a^b G_{\gamma(s)}(\dot{\gamma}(s)) ds,
\end{equation}
among all (piecewise smooth) curves $ \gamma: [a,b] \rightarrow M $ with fixed endpoints (see e.g. \cite{javaloyes2015}). These critical points are preserved under positive affine reparametrizations of $ \gamma $, even if $ F $ is time-dependent. Although we cannot characterize in general pregeodesic curves as critical points of a length functional, as before, since $ \sqrt{G} $ can be non-smooth on lightlike vectors, we can still define the {\em pregeodesics} of $ G $ as the curves that can be parametrized as geodesics (keeping the orientation). Of course, when $ F $ is time-independent, then $ \gamma(s) = (t(s),\sigma(s)) $ is a pregeodesic of $ G $ if and only if $ \sigma(s) $ is a pregeodesic of $ F $ (see \cite[Proposition~B.1]{caponio2018}). This way, we can see pregeodesics of $ G $ (or more specifically, their projection on $ N $) as a generalization of those of $ F $ in the time-dependent case.

\subsection{Traveltime}
\label{subsec:traveltime}
Let us return to the situation we considered in \S~\ref{subsec:speed}, where a wave ray follows a certain trajectory $ \sigma: [a,b] \rightarrow N $ in $ \mathcal{N}_{x,y} $, necessarily satisfying \eqref{eq:time_der}.
\subsubsection{Time-independent case}
Assume first that the speed profile of the wave is time-independent, i.e. $ V = V(x,v) $, so that it provides a proper Finsler metric $ F_x(v) = \frac{||v||_{h_x}}{V(x,v)} $. In this case, note that we can directly integrate \eqref{eq:time_der} to obtain the total traveltime spent by the wave ray when following $ \sigma $:
\begin{equation*}
T[\sigma] \coloneqq t(b)-t(a) = \int_a^b \frac{||\dot{\sigma}(s)||_{h_{\sigma(s)}}}{V(\sigma(s),\dot{\sigma}(s))} ds = \int_a^b F_{\sigma(s)}(\dot{\sigma}(s)) ds = L_F[\sigma],
\end{equation*}
which coincides with the $ F $-length of $ \sigma $, defined in \eqref{eq:F-length}, i.e. the (time-independent) Finsler metric $ F $ acts in a similar way as $ h $ in \S~\ref{subsec:distance}, but instead of measuring distances, it directly measures traveltimes. According to Fermat's principle, $ \sigma $ must be a critical point of $ T $ among all curves in $ \mathcal{N}_{x,y} $, which means that wave rays follow the pregeodesic curves of $ F $.

\subsubsection{Time-dependent case}
When the speed profile is time-dependent, the differential equation \eqref{eq:time_der} is not so easy to solve for the time function $ t(s) $, since now it appears inside $ V $---observe that we cannot integrate \eqref{eq:time_der} directly unless we know $ t(s) $ in advance.\footnote{Although we will not follow this approach here, calculus of variations can be directly applied to \eqref{eq:time_der}, leading to a generalized version of the Euler-Lagrange equations for $ \sigma $ to be a critical point of the traveltime (see e.g. \cite[\S~2]{godin2004} and \cite[\S~3]{levi-civita1931}).} Moreover, as previously discussed, $ F $ is no longer a standard metric whose pregeodesics can be defined in the usual way. Instead, we can consider the more general spacetime viewpoint constructed in \S~\ref{subsec:L}, using the Lorentz-Finsler metric $ G $ on $ M = \mathds{R} \times N $. Note that the curve $ \gamma(s) = (t(s),\sigma(s)) $ in the spacetime directly incorporates the measurement of the traveltime in its first component. Also, we have already pointed out that \eqref{eq:time_der} is equivalent to \eqref{eq:time_der_L}, i.e. $ \gamma $ must be lightlike when representing a wave path. This way, traveling between two points $ x,y \in N $ in the space is equivalent for the wave ray to go from $ (0,x) \in M $ (choosing $ t = 0 $ as the initial time) to the vertical line $ \ell_y: t \mapsto (t,y) $ following a lightlike curve in the spacetime.

Therefore, if we consider the set $ \mathcal{N}_{x,\ell_y} $ of all (piecewise smooth) lightlike curves $ \gamma: [a,b] \rightarrow M $ from $ \gamma(a) = (0,x) $ to $ \gamma(b) \in \ell_y $, the {\em traveltime functional} is given by
\begin{equation}
\label{eq:traveltime_L}
T[\gamma] \coloneqq t(b) = \ell_y^{-1}(\gamma(b))
\end{equation}
and Fermat's principle tells us that a wave path in the spacetime must be a critical point of $ T $ among all curves in $ \mathcal{N}_{x,\ell_y} $. It turns out that these critical points are precisely the lightlike pregeodesics of $ G $ (see \cite[Theorem~4.2]{perlick2006}).

\subsection{Initial wavefront}
\label{subsec:initial_wavefront}
From the previous study, we have come to the conclusion that the path taken by a wave ray between two fixed points is a lightlike pregeodesic of $ G $. Let us assume now, according to the initial formulation of the problem, that the wave departs not from a single point, but from an entire subset $ S \subset N $ that acts as the {\em initial wavefront} at time $ t = 0 $. We will assume that this initial wavefront is a compact (smooth embedded) hypersurface of $ N $ (if $ n = 2 $, then $ S $ is a simple closed curve).\footnote{This can be generalized to the case when the initial wavefront $ S $ is a submanifold of $ N $ of arbitrary codimension (see \cite{pendas2021}).} The wave paths then have to be critical points of the traveltime functional \eqref{eq:traveltime_L}, but now in the broader set $ \mathcal{N}_{S,\ell_y} $ of (piecewise smooth) lightlike curves $ \gamma: [a,b] \rightarrow M $ from $ \sigma(a) \in S $ to $ \gamma(b) \in \ell_y $. This imposes an additional requirement for $ \gamma $ to be critical: $ \gamma $ must depart orthogonally from $ S $ (see \cite[\S~4.1]{pendas2021}).

\begin{defi}
\label{def:orth}
Let $ F $ be a (time-dependent) Finsler metric on $ N $.
\begin{itemize}
\item We say that $ v \in T_xN \setminus \{0\} $ is {\em $ F $-orthogonal} to $ u \in T_xN $ {\em at time} $ t $, denoted $ v \bot_{F_{t,x}} u $, if
\begin{equation*}
g^{F_{t,x}}_v(v,u) = \frac{1}{2} \left. \frac{\partial}{\partial \delta} F_{t,x}(v + \delta u)^2 \right\rvert_{\delta=0} = 0.
\end{equation*}
\item If $ S $ is the initial wavefront (at time $ t = 0 $) and $ x \in S $, we say that $ v \in T_xN \setminus \{0\} $ is {\em $ F $-orthogonal} to $ S $, denoted $ v \bot_F S $, if $ v \bot_{F_{0,x}} u $ for all tangent vectors $ u \in T_xS $.
\end{itemize}
Note that this orthogonality relation is not symmetric, but it is invariant under (positive) vector rescalings: $ v \bot_{F_{t,x}} u $ if and only if $ \lambda_1 v \bot_{F_{t,x}} \lambda_2 u $ for any $ \lambda_1 > 0 $, $ \lambda_2 \in \mathds{R} \setminus \{0\} $.
\end{defi}

\begin{thm}
Let $ \gamma: [a,b] \rightarrow M $ be a curve $ \gamma(s) = (t(s),\sigma(s)) $ departing from $ \sigma(a) \in S $ at time $ t(a) = 0 $. Then, $ \gamma $ is a critical point of the traveltime functional \eqref{eq:traveltime_L} on $ \mathcal{N}_{S,\ell_{\sigma(b)}} $ if and only if it is a lightlike pregeodesic of $ G $ such that $ \dot{\sigma}(a) \bot_F S $.
\end{thm}

\begin{rem}
\label{rem:outward}
Given any point $ x \in S $, there are exactly two $ F $-orthogonal directions to $ S $ (or, equivalently, two $ G $-orthogonal lightlike directions to $ S $; see \cite[Proposition~5.2]{aazami2016}): one points to the interior of $ S $, and the other to the exterior. Here we will focus on the outward trajectories---which are usually the most interesting ones from a practical viewpoint---although the results we will present directly apply to the inward trajectories simply by replacing ``outward'' with ``inward''.
\end{rem}

\subsection{Geodesic equations}
\label{subsec:geodesics}
Now that we have stated the importance of the (pre-)geodesic curves, let us derive the differential equations that will allow us to explicitly compute them. In coordinates, we will use the notation $ p = (t,x) = (x^0,\ldots,x^n) \in M \subset \mathds{R}^{n+1} $, $ \hat{v} = (v^0,\ldots,v^n) \in T_pM \equiv \mathds{R}^{n+1} $ and, consistently,
\begin{equation*}
\begin{split}
G_p(\hat{v}) & = G(x^0,\ldots,x^n,v^0,\ldots,v^n), \\
\gamma(s) & = (\gamma^0(s),\ldots,\gamma^n(s)), \\
\end{split}
\end{equation*}
so that
\begin{equation*}
\dot{\gamma}(s) = (\dot{\gamma}^0(s),\ldots,\dot{\gamma}^n(s)), \quad \text{with} \quad \dot{\gamma}^i(s) \coloneqq \frac{d\gamma^i}{ds}(s) \quad \text{and} \quad \ddot{\gamma}^i(s) \coloneqq \frac{d^2\gamma^i}{ds^2}(s).
\end{equation*}
Also, Einstein's summation convention will be used, summing the indices that appear up and down. Indices $ i,j,k,r $ run from $ 0 $ to $ n $.

As previously defined, a geodesic $ \gamma: [a,b] \rightarrow M $ of $ G $ is a critical point of the $ G $-energy functional \eqref{eq:L-energy}, which means that it satisfies the Euler-Lagrange equations:
\begin{equation}
\label{eq:euler-lagrange}
\frac{d}{ds}\frac{\partial G}{\partial v^r}(\gamma(s),\dot{\gamma}(s)) = \frac{\partial G}{\partial x^r}(\gamma(s),\dot{\gamma}(s)), \quad r=0,\ldots,n.
\end{equation}
Now, observe that
\begin{equation*}
\begin{split}
\frac{\partial G}{\partial v^r}(\gamma(s),\dot{\gamma}(s)) & = 2 \hat{g}_{ri}(\gamma(s),\dot{\gamma}(s))\dot{\gamma}^i(s), \\
\frac{\partial G}{\partial x^r}(\gamma(s),\dot{\gamma}(s)) & = \frac{\partial \hat{g}_{ij}}{\partial x^r}(\gamma(s),\dot{\gamma}(s)) \dot{\gamma}^i(s) \dot{\gamma}^j(s),
\end{split}
\end{equation*}
where $ \hat{g}_{ij}(p,\hat{v}) $ refers to the coefficients of the coordinate matrix of the fundamental tensor $ g^{G_p}_{\hat{v}} $, given by \eqref{eq:gL_coordinates}. Taking the derivative with respect to the parameter $s$ in the first equation, and omitting the explicit reference to $ s $ to avoid cluttering:
\begin{equation*}
\begin{split}
\frac{d}{ds}\frac{\partial G}{\partial v^r}(\gamma,\dot{\gamma}) & = 2\left( \frac{\partial \hat{g}_{ri}}{\partial x^j}(\gamma,\dot{\gamma})\dot{\gamma}^j\dot{\gamma}^i+\hat{g}_{ri}(\gamma,\dot{\gamma})\ddot{\gamma}^i\right) \\ & = 2\left( \frac{1}{2}\frac{\partial \hat{g}_{ri}}{\partial x^j}(\gamma,\dot{\gamma})\dot{\gamma}^j\dot{\gamma}^i + \frac{1}{2}\frac{\partial \hat{g}_{rj}}{\partial x^i}(\gamma,\dot{\gamma})\dot{\gamma}^i\dot{\gamma}^j + \hat{g}_{ri}(\gamma,\dot{\gamma})\ddot{\gamma}^i\right),
\end{split}
\end{equation*}
so, rearranging the terms, \eqref{eq:euler-lagrange} becomes
\begin{equation*}
\hat{g}_{ri}(\gamma,\dot{\gamma}) \ddot{\gamma}^i = -\frac{1}{2}\left( \frac{\partial \hat{g}_{rj}}{\partial x^i}(\gamma,\dot{\gamma}) + \frac{\partial \hat{g}_{ri}}{\partial x^j}(\gamma,\dot{\gamma}) - \frac{\partial \hat{g}_{ij}}{\partial x^r}(\gamma,\dot{\gamma}) \right) \dot{\gamma}^i \dot{\gamma}^j.
\end{equation*}
Finally, multiplying both sides by $ \hat{g}^{kr}(\gamma,\dot{\gamma}) $ (the coefficients of the inverse matrix of $ \hat{g}_{ij}(\gamma,\dot{\gamma}) $) and summing over $r$, we obtain the {\em geodesic equations}:
\begin{equation}
\label{eq:L-geodesics}
\ddot{\gamma}^k(s) = -\hat{\gamma}_{\ ij}^k(\gamma(s),\dot{\gamma}(s))\dot{\gamma}^i(s)\dot{\gamma}^j(s), \quad k = 0,\ldots,n,
\end{equation}
where $ \hat{\gamma}_{\ ij}^k(p,\hat{v}) $ are the {\em formal Christoffel symbols} of $ G $, given by
\begin{equation}
\label{eq:L-christoffel}
\hat{\gamma}_{\ ij}^k(p,\hat{v}) \coloneqq \frac{1}{2} \hat{g}^{kr}(p,\hat{v}) \left( \frac{\partial \hat{g}_{rj}}{\partial x^i}(p,\hat{v}) + \frac{\partial \hat{g}_{ri}}{\partial x^j}(p,\hat{v}) - \frac{\partial \hat{g}_{ij}}{\partial x^r}(p,\hat{v}) \right),
\end{equation}
for any $ p \in M $ and $ \hat{v} \in \mathcal{C}_p $.

\begin{rem}
\label{rem:uniqueness}
Observe that \eqref{eq:L-geodesics} is an ordinary differential equation (ODE) system and thus, geodesics are univocally determined by their initial position and velocity. So, fixing $ x \in S $ and $ v \bot_F S $, there exists a unique lightlike geodesic $ \gamma: [a,b] \rightarrow M $ (assuming $ b $ is maximal), written as $ \gamma(s) = (t(s),\sigma(s)) $, such that $ \gamma(a) = (0,x) $ and $ \dot{\sigma}(a) = v $ (recall Remark~\ref{rem:lightlike}). This curve $ \gamma $ is then a critical point of the traveltime functional \eqref{eq:traveltime_L} on $ \mathcal{N}_{S,\ell_{\sigma(s)}} $, for any $ s \in (a,b] $.
\end{rem}

\section{The model}
\label{sec:model}
We now have all the ingredients needed to answer affirmatively the initial question posed at the beginning of \S~\ref{sec:setting} and compute the entire wave propagation. To this end, we will first construct the theoretical model that provides the evolution of the wavefront over time, and then we will outline how to implement it in practice.

\subsection{The wave trajectories}
\label{subsec:wave_trajectories}
In the previous section, we have identified the paths taken by individual wave rays as lightlike pregeodesics of a Lorentz-Finsler metric $ G $, defined by \eqref{eq:def_L}, departing $ F $-orthogonally from the initial wavefront $ S $. Although we have been working with arbitrary parametrizations so far---in order to show explicitly the importance of the path itself, rather than the specific parametrization we choose---in practice it will be more convenient and natural to parametrize the curves by the time coordinate $ t $. This leads to the following definition.

\begin{defi}
\label{def:wave_trajectory}
Any $ t $-parametrized curve $ \sigma: [0,\tau] \rightarrow N $ such that $ \gamma(t) = (t,\sigma(t)) $ is a lightlike pregeodesic of $ G $, with $ \sigma(0) \in S $ and $ \dot{\sigma}(0) \bot_F S $, will be called a {\em wave trajectory}.
\end{defi}

This way, a wave trajectory $ \sigma(t) $ directly represents the position of a wave ray at time $ t $, and $ \dot{\sigma}(t) $ coincides with actual wave velocity, since \eqref{eq:time_der_F} reduces to
\begin{equation}
\label{eq:F_unitary}
F_{t,\sigma(t)}(\dot{\sigma}(t)) = 1,
\end{equation}
which is equivalent to $ \dot{\sigma}(t) \in \Sigma_{t,\sigma(t)} $. From Remarks~\ref{rem:outward} and \ref{rem:uniqueness}, note that a wave trajectory is uniquely determined by its initial position on the initial wavefront, provided that we only consider outward trajectories. The idea then is to encompass all these trajectories in a single map.

\begin{defi}
\label{def:wavemap}
Given the speed profile of the wave $ V(t,x,v) $ and the initial wavefront $ S $, we define the {\em wavemap} as
\begin{equation}
\label{eq:wavemap}
\begin{array}{rrll}
f \colon & [0,\infty)\times S & \longrightarrow & N \\
& (t,x) & \longmapsto & f(t,x),
\end{array}
\end{equation}
so that $ t \mapsto f(t,x) $ is the unique wave trajectory departing from $ x \in S $ and heading outwards.
\end{defi}

\subsection{The wavefront}
\label{subsec:wavefront}
Until now we have studied individual wave trajectories, but not how they generate the final wavefront. The way the wave behaves to provide this wavefront is governed by {\em Huygens' envelope principle} (see e.g. \cite[\S~2]{markvorsen2016} and \cite[\S~3.3]{pendas2021}):
\begin{quote}
The wavefront $ \mathcal{W}_{\tau} $ at time $ t = \tau $ is given by the envelope of all possible wave trajectories departing from a previous wavefront $ \mathcal{W}_{t_0} $ at time $ t = t_0 $ and running an interval $ \tau-t_0 $.
\end{quote}
In the end, trajectories that remain in the wavefront (i.e., those whose endpoint is an actual point in the wavefront) are the ones that not only are critical points of the traveltime functional but in fact global minima (see \cite[\S~4.1]{pendas2021}). Indeed, if a wave trajectory $ \sigma: [t_0,\tau] \rightarrow N $ does not minimize the traveltime from the wavefront $ \mathcal{W}_{t_0} $ to the endpoint $ \sigma(\tau) $, it means that there exists another wave trajectory from $ \mathcal{W}_{t_0} $ that reaches $ \sigma(\tau) $ in a shorter time. Therefore, $ \sigma(\tau) $ does not belong to the envelope of all possible trajectories departing from $ \mathcal{W}_{t_0} $, which means that it cannot be a point in $ \mathcal{W}_{\tau} $. This motivates the following definition.

\begin{defi}
We say that a wave trajectory $ t \mapsto f(t,x_0) $ is {\em time-minimizing} at time $ t = \tau > 0 $ if, for any fixed $ t_0 \in (0,\tau] $, any other wave trajectory $ t \mapsto f(t,x_1) $ such that $ f(t_0,x_0) = f(t_1,x_1) $ satisfies that $ t_0 \leq t_1 $.\footnote{Essentially, this means that $ t \mapsto (t,f(t,x_0)) $ is a global minimum of the traveltime functional \eqref{eq:traveltime_L} on $ \mathcal{N}_{S,\ell_{f(t_0,x_0)}} $, for any $ t_0 \in (0,\tau] $.}
\end{defi}

Of course, if a trajectory is time-minimizing at $ t = \tau $, then it is so also at any time $ t_0 \in (0,\tau] $. The {\em wavefront} at $ t = \tau $ is then given by
\begin{equation*}
\mathcal{W}_{\tau} = \{ f(\tau,x) \in N: x \in S \text{ and } t \mapsto f(t,x) \text{ is time-minimizing at } t = \tau \},
\end{equation*}
and the following key result ensures that every wave trajectory remains in the wavefront, at least for a small time (see \cite[Theorem~4.8]{pendas2021}).\footnote{In fact, the evolution of the wavefront over (a small) time generates a smooth lightlike hypersurface of the spacetime, as described in \cite{pendas2022}.}

\begin{thm}
\label{thm:cut_points}
There exists a time $ \varepsilon > 0 $ such that every wave trajectory $ t \mapsto f(t,x) $ is time-minimizing at $ t = \varepsilon $.
\end{thm}

This means that the wavefront $ \mathcal{W}_{\tau} $ is directly given by $ S \ni x \mapsto f(\tau,x) $, for any $ \tau \leq \varepsilon $. Beyond this time, however, there might be trajectories in the wavemap that no longer provide points in the wavefront. A point $ x \in N $ beyond which a wave trajectory is no longer time-minimizing is called a {\em cut point} (see \cite[\S~4.1]{pendas2021} and \cite[\S~4.1]{pendas2023a}). These points are usually very interesting from a practical perspective, not only because they allow us to accurately describe the wavefront but also because they mark regions where several wave trajectories tend to converge---e.g., in the case of wildfires, these regions can become extremely dangerous for firefighters, so it is important to detect them beforehand (see the comments in \cite[\S~5.2]{pendas2023a}).

\subsection{The implementation}
\label{subsec:implementation}
Lastly, we will see how this model can be implemented in practice. Assuming we know the speed profile $ V(t,x,v) $---and thus, the (time-dependent) Finsler metric $ F $ given by \eqref{eq:finsler}---and the initial wavefront $ S $, one has to implement the following steps:
\begin{enumerate}
\item Parametrize the initial wavefront and pick a finite number of points there.
\item For each initial point, solve the orthogonality equations to obtain the initial velocity.
\item For each initial point and velocity, solve the pregeodesic equations to obtain the unique wave trajectory running from $ t = 0 $ until $ t = \tau $.
\item For each wave trajectory, check if it is time-minimizing at $ t = \tau $.
\item Obtain the wavefront $ \mathcal{W}_{\tau} $ as an interpolation of the endpoints of time-minimizing wave trajectories. Start again at step 1, using $ \mathcal{W}_{\tau} $ as the initial wavefront.
\end{enumerate}
Let us take a closer look at each step.

\subsubsection{Step 1: The initial wavefront}
In order to deal with the initial wavefront $ S $, first we need to choose a parametrization
\begin{equation*}
\begin{array}{rrll}
\alpha \colon & J \subset \mathds{R}^{n-1} & \longrightarrow & S \subset \mathds{R}^n \\
& \theta = (\theta^1,\ldots,\theta^{n-1}) & \longmapsto & \alpha(\theta) = (\alpha^1(\theta),\ldots,\alpha^n(\theta)),
\end{array}
\end{equation*}
so that the vectors $ \{\frac{\partial \alpha}{\partial \theta^1}(\theta),\ldots,\frac{\partial \alpha}{\partial \theta^{n-1}}(\theta)\} $ form a basis of $ T_{\alpha(\theta)}S $, for all $ \theta \in J $. This parametrization can always be obtained, at least locally---simply by taking local coordinates on $ S $.

\begin{exe}
\label{ex:S}
On $ N = \mathds{R}^2 $, $ S $ must be a simple closed curve, so it always admits a global parametrization of the form $ \alpha: [0,2\pi) \rightarrow S $. In the simplest case, let $ S = \mathds{S}^1 $ be the sphere of radius 1 centered at the origin. Then, we can choose
\begin{equation}
\label{eq:alpha}
\begin{array}{rrll}
\alpha \colon & J = [0,2\pi) & \longrightarrow & S \subset \mathds{R}^2 \\
& \theta & \longmapsto & \alpha(\theta) = (\cos(\theta),\sin(\theta)).
\end{array}
\end{equation}
\end{exe}

Finally, we have to pick a finite number of points $ \{\theta_l\}_{l=1}^{m} \subset J $, so that $ \{\alpha(\theta_l)\}_{l=1}^{m} \subset S $ will be the initial points of the wave trajectories we will compute. Naturally, the wavemap \eqref{eq:wavemap} inherits the parametrization of $ S $ and we can define the {\em discrete wavemap} as
\begin{equation*}
\begin{array}{rrll}
f \colon & [0,\infty)\times \{\theta_l\}_{l=1}^{m} & \longrightarrow & N \\
& (t,\theta_l) & \longmapsto & f(t,\theta_l),
\end{array}
\end{equation*}
where $ t \mapsto f(t,\theta_l) $ is the unique wave trajectory departing from $ \alpha(\theta_l) \in S $ and heading outwards.

\subsubsection{Step 2: The orthogonality equations}
For each fixed $ l \in \{1,\ldots,m\} $, the initial velocity $ v_l \in T_{\alpha(\theta_l)}N $ of the wave trajectory departing from $ \alpha(\theta_l) \in S $ has to be $ F $-orthogonal to $ S $ (recall Definitions~\ref{def:orth} and \ref{def:wave_trajectory}). Since $ \{\frac{\partial \alpha}{\partial \theta^1}(\theta_l),\ldots,\frac{\partial \alpha}{\partial \theta^{n-1}}(\theta_l)\} $ is a basis of $ T_{\alpha(\theta_l)}S $, the vector $ v_l $ has to satisfy the following {\em orthogonality equations}:
\begin{equation}
\label{eq:orth_eq}
g_{v_l}^{F_{0,\alpha(\theta_l)}} \left( v_l,\frac{\partial \alpha}{\partial \theta^i}(\theta_l) \right) = 0, \quad i = 1,\ldots,n-1,
\end{equation}
where the fundamental tensor $ g_v^{F_{t,x}} $ of $ F $ is given by \eqref{eq:gF_coordinates}. This is, of course, not enough to determine the $ n $ components of $ v_l = (v_l^1,\ldots,v_l^n) $. We will also need to fix the length of $ v_l $ and its orientation. Regarding the length, \eqref{eq:F_unitary} tells us that $ v_l $ must be $ F $-unitary:
\begin{equation*}
F_{0,\alpha(\theta_l)}(v_l) = 1.
\end{equation*}
Regarding the orientation, from Remark~\ref{rem:outward} we know that there are always two $ F $-orthogonal directions to $ S $ and we are looking for the one pointing outwards. Therefore, we need to ensure that the determinant
\begin{equation*}
\begin{vmatrix}
v_l^1 & \frac{\partial \alpha^1}{\partial \theta^1}(\theta_l) & \cdots & \frac{\partial \alpha^1}{\partial \theta^{n-1}}(\theta_l) \\
\vdots & \vdots & \vdots & \vdots \\
v_l^n & \frac{\partial \alpha^n}{\partial \theta^1}(\theta_l) & \cdots & \frac{\partial \alpha^n}{\partial \theta^{n-1}}(\theta_l)
\end{vmatrix}
\end{equation*}
is always positive or always negative, depending on the specific parametrization of $ S $ via $ \alpha $.

\begin{exe}
Following Example~\ref{ex:S}, note that $ \alpha $ in \eqref{eq:alpha} is a counterclockwise parametrization of $ S = \mathds{S}^1 $. Hence, in order for $ v_l = (v_l^1,v_l^2) $ to point outwards, we need to check that
\begin{equation*}
\begin{vmatrix}
v_l^1 \ & \ \frac{\partial \alpha^1}{\partial \theta}(\theta_l) \\
v_l^2 \ & \ \frac{\partial \alpha^2}{\partial \theta}(\theta_l)
\end{vmatrix}
=
\begin{vmatrix}
v_l^1 \ & \ -\sin(\theta) \\
v_l^2 \ & \ \cos(\theta)
\end{vmatrix}
> 0.
\end{equation*}
Observe that if we choose a clockwise parametrization, then the determinant must be negative.
\end{exe}

Summing up, we need to solve the equation system \eqref{eq:orth_eq} for $ v_l $, with the restrictions that $ v_l $ must be $ F $-unitary and pointing outwards. This provides a unique initial velocity $ v_l \in T_{\alpha(\theta_l)}N \equiv \mathds{R}^n $.

\subsubsection{Step 3: The pregeodesic equations}
For each fixed $ l \in \{1,\ldots,m\} $, the wave trajectory $ t \mapsto f(t,\theta_l) $ has to be the unique lightlike pregeodesic of $ G $ with initial velocity $ v_l $ (recall Definitions~\ref{def:wave_trajectory} and \ref{def:wavemap}). We have already derived in \S~\ref{subsec:geodesics} the geodesics equations \eqref{eq:L-geodesics} to obtain the unique lightlike geodesic $ \tilde{\gamma}(s) = (t(s),\tilde{\sigma}(s)) $ departing from a fixed initial position and velocity. If we reparametrize $ \tilde{\gamma} $ by the time coordinate, so that
\begin{equation*}
\gamma(t) = \tilde{\gamma} \circ s(t) = (t,\sigma(t)) = (\gamma^0(t),\ldots,\gamma^n(t)) = (t,\sigma^1(t),\ldots,\sigma^n(t)),
\end{equation*}
it is easy to check that the geodesic equations transform into the following $ t $-parametrized {\em pregeodesic equations} (here indices $ i,j,k,r $ run from $ 0 $ to $ n $ and we are omitting the explicit reference to the parameter $ t $ of $ \gamma $); see \cite[\S~4.2]{pendas2021}:
\begin{equation}
\label{eq:pregeog_eq_gamma}
\ddot{\gamma}^k = -\hat{\gamma}_{\ ij}^k(\gamma,\dot{\gamma})\dot{\gamma}^i\dot{\gamma}^j + \hat{\gamma}_{\ ij}^0(\gamma,\dot{\gamma})\dot{\gamma}^i\dot{\gamma}^j\dot{\gamma}^k, \quad k = 0,\ldots,n.
\end{equation}
However, notice that $ \gamma^0(t) = t $ is already determined, so from the definition of $ \hat{\gamma}_{\ ij}^k(p,\hat{v}) $ in \eqref{eq:L-christoffel} and the relationship between the fundamental tensors of $ F $ and $ G $ in \eqref{eq:gL_coordinates}, we can write \eqref{eq:pregeog_eq_gamma} in terms of the spatial curve $ \sigma(t) $ and the (time-dependent) Finsler metric $ F $ (here indices $ \mu,\nu,\xi,\zeta $ run from $ 1 $ to $ n $ and we omit again the explicit reference to the parameter $ t $ of $ \sigma $):
\begin{equation}
\label{eq:pregeog_eq_sigma}
\ddot{\sigma}^{\xi} + g^{\xi\mu}(t,\sigma,\dot{\sigma}) \frac{\partial g_{\mu\nu}}{\partial t}(t,\sigma,\dot{\sigma}) \dot{\sigma}^{\nu} + \dot{\sigma}^{\mu} \dot{\sigma}^{\nu} \left( \gamma^{\xi}_{\ \mu\nu}(t,\sigma,\dot{\sigma}) + \frac{1}{2} \frac{\partial g_{\mu\nu}}{\partial t}(t,\sigma,\dot{\sigma}) \dot{\sigma}^{\xi} \right) = 0,
\end{equation}
where $ g_{\mu\nu}(t,x,v) $ are the coefficients of the coordinate matrix of the fundamental tensor $ g_{v}^{F_{t,x}} $, given by \eqref{eq:gF_coordinates}, $ g^{\mu\nu}(t,x,v) $ are the coefficients of its inverse matrix, and $ \gamma^{\xi}_{\ \mu\nu}(t,x,v) $ are the {\em formal Christoffel symbols} of $ F $, given by
\begin{equation*}
\gamma^{\xi}_{\ \mu\nu}(t,x,v) \coloneqq \frac{1}{2} g^{\xi\zeta}(t,x,v)\left(\frac{\partial g_{\zeta\nu}}{\partial x^{\mu}}(t,x,v)+\frac{\partial g_{\zeta\mu}}{\partial x^{\nu}}(t,x,v)-\frac{\partial g_{\mu\nu}}{\partial x^{\zeta}}(t,x,v)\right),
\end{equation*}
for any $ (t,x) \in M $ and $ v \in T_xN $.

In conclusion, in this step we need to solve the ODE system \eqref{eq:pregeog_eq_sigma} for $ \sigma(t) $ in an interval $ t \in [0,\tau] $, with the initial conditions $ \sigma(0) = \alpha(\theta_l) $ and $ \dot{\sigma}(0) = v_l $. This provides the wave trajectory $ t \mapsto f(t,\theta_l) = \sigma(t) $. Repeating steps 2 and 3 for every initial point $ \{\theta_l\}_{l=1}^m $, we obtain the complete discrete wavemap for $ t \in [0,\tau] $. If computing the individual wave trajectories is enough for our purposes, we can stop here. If, on the other hand, an explicit representation of the wavefront is needed, we have to proceed to the next step.

\subsubsection{Step 4: The cut points}
\label{subsec:step4}
Once we have computed the wavemap until a specific time $ t = \tau $, the set of endpoints $ \{f(\tau,\theta_l)\}_{l=1}^m $ are the ones that will define the wavefront $ \mathcal{W}_{\tau} $. However, as explained in \S~\ref{subsec:wavefront}, we have to consider only the endpoints of the wave trajectories that are time-minimizing at $ t = \tau $, i.e. those that have not reached their cut points before $ \tau $. Cut points can be characterized theoretically (see \cite[Proposition~A.1]{pendas2023a}), but in practice they are very difficult to pinpoint because we can only compute a finite number of wave trajectories (see the discussion in \cite[\S~5.2]{pendas2023a}). Therefore, we need a method to estimate with good precision---the bigger the number of points $ m $, the better the accuracy---whether $ f(\tau,\theta_l) $ actually belongs to $ \mathcal{W}_{\tau} $ or not. The following algorithm provides a good example:
\begin{enumerate}
\item For each fixed $ l_0 \in \{1,\ldots,m\} $, check if the corresponding wave trajectory $ t \mapsto f(t,\theta_{l_0}) $ has intersected a different wave trajectory in the discrete wavemap within the interval $ t \in (0,\tau) $.
\begin{enumerate}
\item If there are no intersections, then $ f(\tau,\theta_{l_0}) \in \mathcal{W}_{\tau} $.
\item Otherwise, let $ \{f(t_i,\theta_{l_0})\}_{i=1}^r $ be the set of all the first intersection points, i.e. $ f(t_i,\theta_{l_0}) $ is the first intersection point between $ t \mapsto f(t,\theta_{l_0}) $ and a different wave trajectory $ t \mapsto f(t,\theta_{l_i}) $, with $ f(t_i,\theta_{l_0}) = f(\tilde{t}_i,\theta_{l_i}) $.
\begin{enumerate}
\item If $ t_i < \tilde{t}_i $ for all $ i \in \{1,\ldots,r\} $, then $ f(\tau,\theta_{l_0}) \in \mathcal{W}_{\tau} $.
\item Otherwise, $ f(\tau,\theta_{l_0}) \notin \mathcal{W}_{\tau} $.
\end{enumerate}
\end{enumerate}
\item Define $ \mathcal{E} \coloneqq \{ l_0 \in \{1,\ldots,m\}: f(\tau,\theta_{l_0}) \in \mathcal{W}_{\tau} \} $.
\end{enumerate}

Essentially, we are considering that the endpoints $ \{f(\tau,\theta_l)\}_{l \in \mathcal{E}} $ belong to the wavefront $ \mathcal{W}_{\tau} $, even though this is an approximation: it is possible that, for some $ l \in \mathcal{E} $, $ t \mapsto f(t,\theta_l) $ is not time-minimizing at $ t = \tau $ because another wave trajectory arrives earlier at $ f(\tau,\theta_l) $, but we have not computed this trajectory within the discrete wavemap.

Nevertheless, the previous algorithm can be difficult to implement and expensive in computing time and power. An alternative option is to take advantage of Theorem~\ref{thm:cut_points} and choose $ \tau $ small enough so that there are no intersections between wave trajectories in the discrete wavemap. Even though this time can be very small in theory, we can always choose carefully the initial points to avoid early intersections (e.g., choosing fewer initial points in regions where the wave trajectories tend to converge), at the cost of losing precision. In this case, $ \mathcal{E} = \{1,\ldots,m\} $.

\subsubsection{Step 5: The final wavefront}
Lastly, the wavefront $ \mathcal{W}_{\tau} $ at time $ t = \tau $ is obtained by interpolating the points $ \{f(\tau,\theta_l)\}_{l \in \mathcal{E}} $ as a compact hypersurface of $ N $. In order to ensure that this interpolation provides a good approximation of the final wavefront, it is usually necessary to implement a previous control check to guarantee that the points $ \{f(\tau,\theta_l)\}_{l \in \mathcal{E}} $ are not too dispersed. In case they are, we can either increase the number of initial points $ m $, or reduce the final time $ \tau $ (or both)---in general, the bigger $ m $ and the smaller $ \tau $, the more precise the wavefront interpolation will be.

Once we have computed $ \mathcal{W}_{\tau} $, we can return to step 1 using $ \mathcal{W}_{\tau} $ as the initial wavefront and $ t = \tau $ as the initial time, in order to find the wavefront at later times.

\subsection{Advantages}
\label{subsec:advantages}
As a concluding remark, we stress the main advantages of the model:
\begin{itemize}
\item {\em Flexibility}: this model is not aimed at a specific type of waves. Rather, it provides a general framework that can be applied not only to determine the (time-, position- and direction-dependent) propagation of classical waves, but also of any physical phenomenon that behaves as such---essentially, that satisfies Fermat's and Huygens' principles. Some paradigmatic examples---where similar models based on Finsler geometry have already been used---are sound waves (see \cite{gibbons2010,gibbons2011}), seismic waves (see \cite{antonelli2003,yajima2009}) and wildfire propagation (see \cite{dehkordi2022,pendas2023a,pendas2024,markvorsen2016}). Furthermore, this model can serve as an initial template upon which additional effects can be easily added, such as Snell's law of refraction (see \cite{markvorsen2023}).

\item {\em Efficiency}: the spirit of the model is to obtain the evolution of the wave by computing individual wave trajectories. This simplifies the direct approach of the classical wave equation and its generalizations---e.g., Richards' equations for an elliptical propagation (see \cite{richards1990})---where the wavefront is directly characterized and computed by means of a partial differential equation (PDE) system. Therefore, this model essentially reduces the general PDE system provided by Huygens' principle (see \cite[Theorem~4.14]{pendas2021}) to the computation of several ODE ones---the geodesic equations---which are easier to solve and more efficient, computationally speaking. Moreover, this approach also simplifies the management of cut points and possible crossovers of the wave: the PDE system cannot be solved beyond cut points (see \cite[Remark~4.15]{pendas2021}) and the wavefront has to be corrected each time one appears, whereas the individual wave trajectories can all be computed until the selected final time, regardless of the number of intermediate cut points, and then a single algorithm such as the one described in \S~\ref{subsec:step4} can be implemented to correct the wavefront (see also the discussion in \cite[Appendix~A.1]{pendas2023a}).

\item {\em Accuracy}: the model can handle any propagation pattern and any time, position and direction dependence, so long as the infinitesimal spread (i.e., the indicatrix $ \Sigma_{t,x} $) remains strongly convex. This is a major improvement over models that need to impose some restrictions on the wave spread, either because the differential equations they solve are only valid for a specific propagation shape---e.g., current fire growth simulators based on Richards' equations need to assume that the infinitesimal fire spread always takes the form of an ellipse (see e.g. \cite{finney1998,tymstra2010})---or because this simplification is crucial for the model to be computationally efficient. In our model, every modification of the wave propagation is encoded within the (time-dependent) Finsler metric $ F $ and the geodesic equations (along with the method used to solve them) are formally the same, regardless of how intricate $ F $ is---i.e., increasing the complexity of $ F $ does not entail a technical complication nor does it make the model inefficient.
\end{itemize}

\newpage

\section*{Acknowledgments}
The author was partially supported by:
\begin{itemize}
\item Project PID2021-124157NB-I00, funded by MCIN/AEI/10.13039/5011000 11033/ and ``ERDF A way of making Europe".
\item Project PID2020-116126GB-I00, funded by MCIN/AEI/10.13039/5011000 11033/.
\item The framework IMAG-Mar\'{i}a de Maeztu grant CEX2020-001105-M/AEI/ 10.13039/501100011033.
\item Ayudas a proyectos para el desarrollo de investigaci\'{o}n cient\'{i}fica y t\'{e}cnica por grupos competitivos (Comunidad Aut\'{o}noma de la Regi\'{o}n de Murcia), included in the Programa Regional de Fomento de la Investigaci\'{o}n Cient\'{i}fica y T\'{e}cnica (Plan de Actuaci\'{o}n 2022) of the Fundaci\'{o}n S\'{e}neca-Agencia de Ciencia y Tecnolog\'{i}a de la Regi\'{o}n de Murcia, REF. 21899/PI/22.
\item Ayudas para la Formaci\'{o}n de Profesorado Universitario (FPU) from the Spanish Government.
\end{itemize}

%
%

\end{document}